\documentclass[12pt]{amsart}

\usepackage{amsmath,amscd}
\usepackage{amsthm}
\usepackage{amssymb}
\usepackage{amsfonts}
\usepackage{mathrsfs}
\usepackage{latexsym}
\usepackage{graphicx}
\usepackage{epsfig}
\usepackage[T1]{fontenc}
\usepackage{times}

\theoremstyle{plain}
\newtheorem{theorem}{Theorem}[section]
\newtheorem{lemma}[theorem]{Lemma}
\newtheorem{proposition}[theorem]{Proposition}
\newtheorem{corollary}[theorem]{Corollary}

\theoremstyle{definition}

\newtheorem*{remark}{Remark}

\allowdisplaybreaks

\def\Xint#1{\mathchoice
  {\XXint\displaystyle\textstyle{#1}}%
  {\XXint\textstyle\scriptstyle{#1}}%
  {\XXint\scriptstyle\scriptscriptstyle{#1}}%
  {\XXint\scriptscriptstyle\scriptscriptstyle{#1}}%
  \!\int}
\def\XXint#1#2#3{{\setbox0=\hbox{$#1{#2#3}{\int}$}
    \vcenter{\hbox{$#2#3$}}\kern-.5\wd0}}

\def\dashint{\Xint-}

\renewcommand{\leq}{\leqslant}
\renewcommand{\geq}{\geqslant}

\newcommand{\inner}[2]{\langle #1\,,#2\rangle}
\newcommand{\biginner}[2]{\left\langle #1\,,\,#2\right\rangle}

\newcommand{\ddl}[2]{\frac{d{#1}}{d{#2}}}
\newcommand{\ppl}[2]{\frac{\partial{#1}}{\partial{#2}}}
\newcommand{\ppz}[2]{\frac{\partial^2{#1}}{\partial{#2}^2}}

\newcommand{\C}{\mathbb{C}}
\newcommand{\R}{\mathbb{R}}
\renewcommand{\H}{\mathbb{H}}

\newcommand{\Fcal}{\mathcal{F}}

\newcommand{\Lcal}{\mathcal{L}}

\newcommand{\Xfrak}{\mathfrak{X}}

\newcommand{\I}{\mathbf{I}}

\DeclareMathOperator{\Area}{Area}

\DeclareMathOperator{\Isom}{Isom}
\DeclareMathOperator{\Ric}{Ric}
\DeclareMathOperator{\Rm}{Rm}

\DeclareMathOperator{\dist}{dist}
\DeclareMathOperator{\tr}{tr}

\renewcommand{\div}{\mathop{\mathrm{div}}}
\renewcommand{\mod}{\mathop{\mathrm{mod}}}

\newcommand{\gbar}{\bar{g}}

\newcommand{\Rbar}{{\overline{R}}}

\newcommand{\Ricbar}{{\overline{\Ric}}}

\newcommand{\nablabar}{{\overline{\nabla}}}

\begin{document}

\title{Foliations For Quasi-Fuchsian $3$-Manifolds}

\author{Biao Wang}
\date{\today}
\thanks{Research is partially supported by the National
        Science Foundation grant DMS-0513436}
\subjclass[2000] {Primary 53C44, Secondary 57M05}
\address{Department of Mathematics,
         Cornell University, Ithaca, NY 14853}
\email{wang@math.cornell.edu}

\begin{abstract}In this paper, we prove that if a
quasi-Fuchsian $3$-manifold contains a minimal surface whose
principle curvature is less than $1$, then it admits a foliation
such that each leaf is a surface of constant mean curvature. The
key method that we use here is volume preserving mean 
curvature flow.
\end{abstract}

\maketitle


\section{Introduction}

A codimension one foliation $\Fcal$ of a Riemanian manifold is
called a {\em CMC foliation}, if each leaf of the foliation is a
hypersurface of constant mean curvature. A quasi-Fuchsian group
$\Gamma$ is a Kleinian group which is obtained by a quasiconformal
deformation a Fuchsian group, its limit set is a closed Jordan
curve dividing the domain of discontinuity $\Omega$ on 
$S_{\infty}^{2}$ into two simply connected, invariant component.
Topologically, $(\H^{3}\cup\Omega)/\Gamma=S\times[0,1]$,
where $S$ is a closed surface with $\pi_{1}(\Sigma)=\Gamma$.
In this paper, we always assume that $S$ is a {\em closed}
Riemann surface with genus $\geq{}2$.

Suppose $M$ is a $3$-dimensional quasi-Fuchsian hyperbolic 
manifold, Mazzeo and Pacard proved that each end of $M$ admits a
unique CMC foliation (cf. \cite{MP2007}). Next we may ask if
the whole quasi-Fuchsian manifold $M$ admits a CMC foliation?
If $M$ admits a CMC foliation $\Fcal$, then the foliation $\Fcal$
must contain a leaf $L$ whose mean curvature is zero, i.e. $L$
is a minimal surface in $M$. Therefore we need to know whether
$M$ contains a minimal surface at first. There are several ways 
to prove that $M$ contains a least area minimial surface 
$\Sigma$ with $\pi_{1}(M)\cong\pi_{1}(\Sigma)$ 
(cf. \cite{A1983,MSY1982,SY1979,U1983}). 

In this paper, we will prove the following theorem.

\begin{theorem}\label{thm:main theorem}
Suppose that $M$ is a quasi-Fuchsian $3$-manifold, which contains 
a closed immersed minimal surface $\Sigma$ with genus $\geq{}2$ 
such that $\pi_{1}(M)\cong\pi_{1}(\Sigma)$, if the principle 
curvature $\lambda$ of $\Sigma$ satisfies $|\lambda(x)|<1$ for all 
$x\in\Sigma$, then $M$ admits a unique CMC foliation.
\end{theorem}

We will use the volume preserving mean curvature flow developed
by G. Huisken (cf. \cite{H1984,H1987}) to prove Theorem
\ref{thm:main theorem} in $\S${}\ref{sec:foliation}. This idea 
is inspired by Ecker and Huisken's paper \cite{EH1991}. 
Furthermore, we will show that $M$ doesn't admit a CMC foliation if
the principle curvature of $\Sigma$ is very large in
$\S${}\ref{sec:counterexam}, where the idea
of using infinite minimal catenoids as barrier surfaces 
contributes to Bill Thurston.

This paper is organized as follows. 
In $\S${}\ref{sec:preliminary}, we give some
definitions and basic properties about quasi-Fuchsian groups 
and submanifolds. 
In $\S${}\ref{sec:MCF}, we discuss the volume preserving mean
curvature flow and prove the existence of the long time solution. 
In $\S${}\ref{sec:foliation}, we will prove 
Theorem \ref{thm:main theorem}. 
In $\S${}\ref{sec:counterexam}, we will give a counterexample. 

\subsection*{Acknowledgements} This paper is supervised under Bill 
Thurston. I am grateful to him for his guidance and a lot of helpful 
and stimulating conversations. I also appreciate John Hubbard 
and Xiaodong Cao, who give me many suggestions.
 
\section{Preliminaries}\label{sec:preliminary}

In this section, we review some basic facts on 
quasi-Fuchsian $3$-manifolds and geometry of submanifolds.

\subsection{Quasifuchsian groups} 
A subgroup $\Gamma$ of $\Isom(\H^{3})$ is called a 
{\em Kleinian groups} if $\Gamma$ acts on $\H^{3}$ properly
discontinuously. For any Kleinian group $\Gamma$, 
$\forall\,p\in\H^{3}$, the orbit set
\begin{equation*}
   \Gamma(p)=\{\gamma(p)\ |\ \gamma\in{}\Gamma\}
\end{equation*}
has accumulation points on $S^{2}_{\infty}=\partial\H^{3}$,
these points are called the {\em limit points} of $\Gamma$, 
and the closed set of all these points is called the 
{\em limit set} of $\Gamma$,
which is denoted by $\Lambda_{\Gamma}$. The complement 
of the limit set, i.e.,
\begin{equation*}
   \Omega_{\Gamma}=S^{2}_{\infty}\setminus\Lambda_{\Gamma}\ ,
\end{equation*}
is called the {\em region of discontinuity}. If 
$\Omega_{\Gamma}=\emptyset$, $\Gamma$ is called a Kleinian 
group of the first kind, and otherwise of the second kind.

Suppose $\Gamma$ is a finitely generated torsion free Kleinian 
group which has more than two limit points, we call $\Gamma$
{\em quasi-Fuchsian} if its limit set $\Lambda_{\Gamma}$ is a
closed Jordan curve and both components $\Omega_{1}$ and 
$\Omega_{2}$ of its region of discontinuity are invariant under
$\Gamma$. The limit set $\Lambda_{\Gamma}$ of the quasi-Fuchsian
group $\Gamma$ is either a (standard) circle or a closed Jordan 
curve which fails to have a tangent on an everywhere dense set 
(cf. \cite[Theorem 4.2]{Lehto1987}. When $\Lambda_{\Gamma}$ 
is a circle, we call $\Gamma$ a Fuchsian group. Of course, 
$\Lambda_{\Gamma}$ is invariant under $\Gamma$ too. The following 
statement about quasi-Fuchsian groups can be found in 
\cite[page 8]{CEG2006}.

\begin{proposition}[{\bf Maskit}\ \cite{Maskit1970}, 
{\bf Thurston}\ \cite{Thurston1980}]
If $\Gamma$ is a finitely generated, torsion-free Kleinian group,
then the following conditions are equivalent:
\begin{enumerate}
   \item $\Gamma$ is quasi-Fuchsian.
   \item $\Omega_{\Gamma}$ has exactly two components, each of
         which is invariant under $\Gamma$.
   \item There exist a Fuchsian group $G$ and a 
         quasiconformal homeomorphism 
         $w:\widehat{\C}\to\widehat{\C}$ such that
         $\Gamma=w\circ{}G\circ{}w^{-1}$.
\end{enumerate}
\end{proposition}

For a finitely generated, torsion free
quasi-Fuchsian group $\Gamma$ with invariant components 
$\Omega_{1}$, $\Omega_{2}$ of $\Omega_{\Gamma}$, Albert Marden
(cf. \cite{Marden1974}) proved that $\Gamma$ has the following
properties:
\begin{itemize}
   \item Each of $S_{1}=\Omega_{1}/\Gamma$ and  
         $S_{2}=\Omega_{2}/\Gamma$ is a finitely 
         punctured Riemann surface.
   \item $M_{\Gamma}=\H^{3}/\Gamma$ is diffeomorphic to
         $(\Omega_{1}/\Gamma)\times(0,1)$, and
         $\overline{M}_{\Gamma}=(\H^{3}
         \cup\Omega_{\Gamma})/\Gamma$ is 
         diffeomorphic to $(\Omega_{1}/\Gamma)\times[0,1]$.
\end{itemize}
We will call $M_{\Gamma}$ a {\em quasi-Fuchsian $3$-manifold}. 
In this paper we write $M_{\Gamma}=S\times\R$, where $S$
is a closed surface with genus $\geq{}2$.

\subsection{Geomerty of submanifolds} 
In this subsection, we rephrase some materials from 
\cite{U1983} for convenience. Let $(M,\gbar_{\alpha\beta})$ be a
quasi-Fuchsian $3$-manifold, and let $\Sigma$ be a 
immersed minimal surface in $M$. Suppose the coordinate system on
$\Sigma\equiv\Sigma\times\{0\}$ is isothermal so that the induced
metric $g=(g_{ij})_{2\times{}2}$ on $\Sigma$ can be written in
the form
\begin{equation*}
   g(x,0)=\{g_{ij}(x,0)\}_{1\leq{}i,j\leq{}2}=
   e^{2v(x)}\I
\end{equation*}
where $\I$ is a $2\times{}2$ unit matrix, and let
\begin{equation*}
   A(x)\equiv{}A(x,0)=\{h_{ij}(x,0)\}
\end{equation*}
be the second fundamental form of $\Sigma$.

In a collar neighborhood of $\Sigma$ in $M$, there exists normal
coordinates induced by $\exp:T^{\bot}\Sigma\to{}M$ in a neighborhood
on which
\begin{equation*}
    \Sigma\times(-\varepsilon,\varepsilon)
    \subset{}T^{\bot}\Sigma\to{}M
\end{equation*}
is a (local) diffeomorphism. If coordinates 
$(x^{1},x^{2})$ are introduced on $\Sigma$, then
\begin{equation*}
   \exp((x^{1},x^{2}),x^{3})=(x^{1},x^{2},x^{3})
\end{equation*}
induces a coordinate patch in $M$. Choose
$p=(x^{1},x^{2},x^3)=(x,r)$ the local coordinate system in a
neighborhood of $\Sigma$ so that $\Sigma=\{(x,r)\in{}M\ |\ r=0\}$.
Let $N_{0}$ be the unit normal vector field on $\Sigma$, and let
\begin{equation}\label{eq:Surface with distance r}
   \Sigma(r)=\{\exp_{x}r{}N_{0}\ |\ x\in\Sigma\}
\end{equation}
for a small positive constant $r$. For
$(x,r)\in\Sigma\times(-\varepsilon,\varepsilon)
\subset{}T^{\bot}\Sigma$,
it's well known that the pullback metric has the form
\begin{equation}\label{eq:Metric of Product}
   \gbar(x,r)=\begin{pmatrix}
                 g(x,r) & 0\\
                      0 & 1
              \end{pmatrix}
             =\begin{pmatrix}
                 g_{11}(x,r) & g_{12}(x,r) & 0 \\
                 g_{21}(x,r) & g_{22}(x,r) & 0 \\
                           0 &           0 & 1
              \end{pmatrix}
\end{equation}
where $g(x,r)$ is the induced metric on $\Sigma(r)$.

The second fundamental form $A=(h_{ij})$ of
$\Sigma(r)$ is a $2\times{}2$ matrix defined by
   \begin{equation}
      h_{ij}=\inner{\nablabar_{e_{i}}e_{3}}{e_{j}}\ ,
      \quad{}1\leq{}i,j\leq{}2\ ,
   \end{equation}
where $\nablabar$ is the covariant differentiation in $M$, and
$\{e_{1},e_{2},e_{3}\}$ is the local frame for $M$ such that $e_{3}$
is the unit normal vector of $\Sigma(r)$ and $e_{1},e_{2}$ are two
unit vectors in the tangent plane of $\Sigma(r)$. Direct computation
shows that the second fundamental forms $A(x,r)=\{h_{ij}(x,r)\}$ on
$\Sigma(r)$ are given by
\begin{equation}\label{eq:second fundamental form}
   h_{ij}(x,r)=\frac{1}{2}\,\ppl{}{r}\,g_{ij}(x,r)\ ,
   \quad{}1\leq{}i,j\leq{}2\ .
\end{equation}
Note that the sectional curvature of $M$ is $-1$, there are three
curvature equations of the form
\begin{equation}\label{eq:curvature operator}
   \Rbar_{i3j3}=-(\gbar_{33}\gbar_{ij}-\gbar_{i3}\gbar_{3j})
   =-g_{ij}\ ,\quad{}1\leq{}i,j\leq{}2\ ,
\end{equation}
where the Riemann curvature tensor is given by
\begin{equation*}
   \Rbar(X,Y)Z=-\nablabar_{X}\nablabar_{Y}Z+
                \nablabar_{Y}\nablabar_{X}Z+
                \nablabar_{[X,Y]}Z
\end{equation*}
for $X,Y,Z\in\Xfrak(M)$. Direct computation shows that the curvature
forms are given by
\begin{equation}\label{eq:curvature equation}
      \Rbar_{i3j3}=\frac{1}{2}\,\ppz{g_{ij}}{r}-
      \frac{1}{4}\,g^{kl}\,\ppl{g_{il}}{r}
      \ppl{g_{jk}}{r}\ ,\quad{}1\leq{}i,j\leq{}2\ .
\end{equation}
From \eqref{eq:curvature operator} and
\eqref{eq:curvature equation}, we get partial differential
equations
\begin{equation}
   -g_{ij}=\frac{1}{2}\,\ppz{g_{ij}}{r}-
  \frac{1}{4}\,g^{kl}\,\ppl{g_{il}}{r}
  \ppl{g_{jk}}{r}\ ,
\end{equation}
whose solutions can be written in the form
\begin{equation}\label{eq:Metric of QF}
   g(x,r)=e^{2v(x)}[\cosh{}r{}\,\I+
   \sinh{}r{}e^{-2v(x)}A(x)]^{2}
\end{equation}
for all $x\equiv(x,0)\in\Sigma$ and $-\epsilon<r<\epsilon$.
This metric is nonsingular in a collar neighborhood of $\Sigma$ in
any case. If the principle curvature of $\Sigma\subset{}M$
\begin{equation*}
   \lambda(x)=\sqrt{-\det\,[A(x)e^{-2v(x)}]}<1\ ,
\end{equation*}
then it is non-singular for all $r\in\R$. 

\begin{proposition}The mean curvature of $\Sigma(r)$ is given by
\begin{equation}\label{eq:mean curvature of Sigma(r)}
   H(x,r)=\frac{2(1-\lambda^{2}(x))\tanh{}r}
   {1-\lambda^{2}(x)\tanh^{2}r}\ ,
   \quad\forall\,x\in\Sigma\ ,
\end{equation}
here the normal vector on $\Sigma(r)$ points to the minimal 
surface $\Sigma$.
\end{proposition}

\begin{proof}In order to compute the mean curvature $H$, we
need to find the eigenvalues of the second fundamental form
$A(x,r)$. In other words, we need solve the equation
\begin{equation*}
   \det\,[h_{ij}-\mu{}g_{ij}]=0\ ,
\end{equation*}
which is equivalent to the equation
\begin{equation*}
   \det\,[(\sinh{}r\I+\cosh{}r{}e^{-2v(x)}A(x))-
   \mu(\cosh{}r\I+\sinh{}r{}e^{-2v(x)}A(x))]=0\ .
\end{equation*}
Solve the above equation, we get two eigenvalues:
\begin{equation*}
   \mu_{1}=\frac{\tanh{}r-\lambda(x)}{1-\lambda(x)\tanh{}r}
   \qquad\text{and}\qquad
   \mu_{2}=\frac{\tanh{}r+\lambda(x)}{1+\lambda(x)\tanh{}r}\ .
\end{equation*}
Since $H=\mu_{1}+\mu_{2}$, the proposition follows.
\end{proof}

It's easy to check that $H(x,r)$ defined in
\eqref{eq:mean curvature of Sigma(r)} is a monotonically 
increasing function with respect to $r$, i.e. 
$H(x,r_{1})\leq{}H(x,r_{2})$ if $r_{1}\leq{}r_{2}$. In fact, 
we have
\begin{equation*}
   \ppl{}{r}\,H(x,r)=
   \frac{2(1-\lambda^{2}(x))[1+\lambda^{2}(x)\tanh^{2}r]}
   {[1-\lambda^{2}(x)\tanh^{2}r]^{2}\cosh^{2}r}\geq{}0\ ,
   \quad\forall\,x\in\Sigma\ .
\end{equation*}
As $r\to\pm\infty$, $H\to\pm{}2$, and as $r\to{}0$, $H\to{}0$.

\begin{theorem}[{\bf Uhlenbeck}\ \cite{U1983}]
\label{thm:U1983-(3.3)}
If $M$ is a complete, hyperbolic manifold and $\Sigma$ is a minimal
surface in $M$ with $|\lambda(x)|<1$ for all $x\in{}\Sigma$, then
   \begin{enumerate}
      \item $\exp{}T^{\bot}\Sigma\cong\widetilde{M}\to{}M$,
            where $\widetilde{M}$ is the cover of $M$
            corresponding to
            $\pi_{1}(\Sigma)\subset\pi_{1}(M)$.
      \item $\widetilde{M}$ is quasi-Fuchsian.
      \item $\Sigma\subset{}M$ is area minimizing;
            $\Sigma\subset\widetilde{M}$ is the
            only closed minimal surface of any type in
            $\widetilde{M}$.
      \item $\Sigma\subset\widetilde{M}$ is embedded.
      \item $\Sigma\subset{}M$ is totally geodesic if and only
            if $\widetilde{M}$ is Fuchsian.
   \end{enumerate}
\end{theorem}

\begin{corollary}\label{cor:corollary of Uhlenbeck}
Suppose $\Sigma$ is an immersed minimal
surface in a quasi-Fuchsian $3$-manifold $M$ which is
homotopic to $\Sigma$, if the principle curvature of $\Sigma$ is
between $-1$ and $1$, then
\begin{itemize}
   \item $\Sigma$ is the unique minimal surface which is
         embedded in $M$, 
   \item the metric $\gbar_{\alpha\beta}$ on $M=\Sigma\times\R$ 
         is given by \eqref{eq:Metric of Product} and 
         \eqref{eq:Metric of QF}, and
   \item $M$ can be foliated by either the geodesics perpendicular
         to the minimal surface $\Sigma$ or the equidistant 
         surfaces $\{\Sigma(r)\}_{-\infty<r<\infty}$ defined by
         \eqref{eq:Surface with distance r}.
\end{itemize}
\end{corollary}

\section{Volume preserving mean curvature flow}\label{sec:MCF}

In this section, we will discuss the volume preserving mean 
curvature flow developed by G. Huisken and others. A good 
reference for mean curvature flow is the book written by 
Xi-Ping Zhu (cf. \cite{Zhu2002}).

By the discussion in $\S${}\ref{sec:preliminary},
$(M,\gbar_{\alpha\beta})$ can be foliated either by the geodesics
which are perpendicular to the minimal surface $\Sigma$ or by the
surfaces $\Sigma(r)$ for all $r\in\R$, where $\Sigma(r)$ is defined
by \eqref{eq:Surface with distance r}. Denote by $N$ the unit
tangent vector field on the geodesics, which is a well defined
vector field on $M$.

For any tensor field $\Phi$ on $(M,\gbar_{\alpha\beta})$ we define
the supremum norms by
\begin{equation*}
   \|\Phi\|=\sup_{x\in{}M}|\Phi(x)|_{\gbar_{\alpha\beta}}
   \qquad\text{and}\qquad
   \|\Phi\|_{k}=\sum_{j=0}^{k}\|\nablabar{}^{j}\Phi\|\ .
\end{equation*}

\subsection{Evolution equations}
Let $S$ be a smooth surface which is diffeomorphic to the minimal
surface $\Sigma\subset{}M$, and let $F_{0}^{r}:S\to{}M$ be 
the immersion of $S$ in $M$ such that $F_{0}^{r}(S)=\Sigma(r)$ 
for some positive constant $r$. Next we consider a family of 
smoothly immersed surfaces in $M$,
\begin{equation*}
   F:S\times[0,T)\to{}M\ ,
   \quad{}0\leq{}T\leq\infty
\end{equation*}
with $F(\cdot,0)=F_{0}^{r}$. For each $t\in[0,T)$, write
\begin{equation*}
   S_{t}=S_{t}(r)=\{F(x,t)\in{}M\ |\ x\in{}S\}\ .
\end{equation*}
We need define some quantities and operators on $S_{t}$:
\begin{itemize}
   \item the induced metric of $S_{t}$ is denoted by
         $g=\{g_{ij}\}$,
   \item the second fundamental form of $S_{t}$ is
         denoted by $A=\{h_{ij}\}$,
   \item the mean curvature of $S_{t}$ with respect to the
         normal pointing to the minimal surface $\Sigma$ is given 
         by $H=g^{ij}h_{ij}$,
   \item the square norm of the second fundamental form of
         $S_{t}$ is given by
         \begin{equation*}
            |A|^{2}=g^{ij}g^{kl}h_{ik}h_{jl}\ ,
         \end{equation*}
   \item the covariant derivative of $S_{t}$ is denoted by
         $\nabla$,
   \item the Laplacian on $S_{t}$ is given by
         $\Delta=g^{ij}\nabla_{i}\nabla_{j}$.
\end{itemize}
Each quantity or operator with respect to
$(M,\gbar_{\alpha\beta})$ will be added a bar on its top. The 
curvature operator $\overline{\Rm}$ on $(M,\gbar_{\alpha\beta})$ 
is given by
\begin{equation}
   \Rbar_{\alpha\beta\gamma\delta}=
   -(\gbar_{\alpha\gamma}\gbar_{\beta\delta}-
   \gbar_{\alpha\delta}\gbar_{\beta\gamma})\ ,
   \quad{}1\leq\alpha,\beta,\gamma,\delta\leq{}3\ .
\end{equation}

We consider the volume preserving mean curvature 
flow (cf. \cite{H1987}):
\begin{equation}\label{eq:Vol-MCF}
   \left\{
   \begin{aligned}
      &\ppl{}{t}\,F(x,t)=[h(t)-H(x,t)]\nu(x,t)\ ,
      \quad{}x\in{}S\ ,\ 0\leq{}t<T\ ,\\
      &F(\cdot,0)=F_{0}^{r}\ ,
   \end{aligned}
   \right.
\end{equation}
where
\begin{equation*}
   h(t)=\dashint_{S_{t}}H{}d\mu=
   \frac{1}{\Area(S_{t})}\int_{S_{t}}H{}d\mu
\end{equation*}
is the {\em average mean curvature} of $S_{t}$, and 
$\nu$ is the normal on $S_{t}$ so that $-\nu$ points to the minimal surface $\Sigma$. It's easy to verify that 
the volume of the domain bounded by $\Sigma$ and $S_{t}$ is
independent of time. In \cite{H1986,H1987}, Huisken proved the
following theorem.

\begin{theorem}[{\bf Huisken}]\label{thm:Huisken1986-(7.1)}
If the initial surface $S_{0}$ is smooth, then \eqref{eq:Vol-MCF}
has a smooth solution on some maximal open time interval
$0\leq{}t<T$, where $0<T\leq\infty$. If $T<\infty$, then
\begin{equation}
   |A|_{\max}(t)\equiv\max_{x\in{}S}|A|(x,t)\to\infty\ ,
   \quad\text{as}\ t\to{}T\ .
\end{equation}
\end{theorem}

In this section, we will prove the following theorem.

\begin{theorem}\label{thm:long time existence}
For any fixed $r>0$, the evolution equation \eqref{eq:Vol-MCF} 
has a unique long time solution $(i.e. T=\infty)$. 
As $t\to\infty$, the surfaces $\{S_{t}\}$ converge exponentially
fast to a smooth surface $S_{\infty}$ of constant mean curvature.
\end{theorem} 

For this aim, we assume $T<\infty$ at the very beginning, if 
we can prove that there exist constants 
$\{C(m)\}_{m=0,1,2,\ldots}$ independent of time such that the 
estimates
\begin{equation}
   |\nabla^{m}A|^{2}\leq{}C(m)\ , \quad{}m=0,1,2,\ldots\ ,
\end{equation}
are uniformly on $S_{t}$ for $0\leq{}t<T$, then we can derive 
that the limit surface
$S_{T}=\displaystyle\lim_{t\to{}T}S_{t}$ is a smooth surface, so
we can extend $T$ a little bit further by Theorem
\ref{thm:Huisken1986-(7.1)}, this is contradicted to the hypothesis
that $T$ is maximal.

To obtain in the next step a priori estimate for $|A|^{2}$, we need
evolution equations for the metric and the second fundamental form
on $S_{t}$.

\begin{lemma}[{\bf Huisken--Yau}\ \cite{HY1996}]
\label{lem:HY1996-(3.6)}
We have the following evolution equations:
\begin{enumerate}
   \item $\displaystyle\ppl{}{t}\,g_{ij}=2(h-H)h_{ij}$,
   \item $\displaystyle\ppl{}{t}\,h_{ij}=\nabla_{i}\nabla_{j}H+
         (h-H)h_{il}g^{kl}h_{kj}+(h-H)g_{ij}$,
   \item $\displaystyle\ppl{}{t}\,\nu=\nabla{}H$,
   \item $\displaystyle\ppl{}{t}\,\mu=H(h-H)\mu$, where $\mu$ is
         the measure on $S_{t}$.
\end{enumerate}
\end{lemma}

Since $(M,\gbar_{\alpha\beta})$ is a $3$-manifold with
constant sectional curvature, we have 
$\nablabar_{m}\Rbar_{ijkl}\equiv{}0$, 
$\Ricbar(\nu,\nu)=-2$, and
\begin{equation*}
   h_{ij}h_{jl}\Rbar_{lmlm}-h_{ij}h_{lm}\Rbar_{limj}=
   -(\lambda_{1}-\lambda_{2})^{2}=H^{2}-2|A|^{2}\ .
\end{equation*}
Together with  Simons' identity (cf. 
\cite[Lemma 1.3(i)]{HY1996}), we 
obtain the following additional evolution equations.

\begin{lemma}[{\bf Huisken--Yau}\ \cite{HY1996}]
\label{lem:HY1996-(3.7)}
Under the evolution equation \eqref{eq:Vol-MCF}, the second
fundamental form satisfies the evolution equations
\begin{enumerate}
   \item $\displaystyle\ppl{}{t}\,h_{ij}=\Delta{}h_{ij}+
         (h-2H)h_{il}g^{lk}h_{kj}+(|A|^{2}+2)h_{ij}+
         (h-2H)g_{ij}$,
   \item $\displaystyle\ppl{}{t}\,H=\Delta{}H+
         (H-h)(|A|^{2}-2)$,
   \item $\displaystyle\ppl{}{t}\,|A|^{2}=\Delta|A|^{2}-
          2|\nabla{}A|^{2}+2|A|^{4}-2h\tr{}A^{3}+
          4|A|^{2}+2H(h-2H)$, where 
          $\tr{}A^{3}=\dfrac{H}{2}(3|A|^{2}-H^{2})$.
\end{enumerate}
\end{lemma}

\subsection{Existence of the long time solution} 
Define a function $\ell:M\to\R$ by
\begin{equation*}
   \ell(p)=\dist(p,\Sigma)=\min\{\dist(p,p')\ |\ p'\in\Sigma\}
\end{equation*}
for all $p\in{}M$, where $\dist(\cdot,\cdot)$ is the distance
function on $(M,\gbar_{\alpha\beta})$. By Corollary
\ref{cor:corollary of Uhlenbeck}, every point $p\in{}M$ has 
the form $p=(p',r)$ for some point $p'\in\Sigma$, where 
$r=\ell(p)$. Let
\begin{equation*}
   u=\ell|{S_{t}}
   \qquad\text{and}\qquad
   \Theta=\inner{N|{S_{t}}}{\nu}
\end{equation*}
be the height function and the gradient function of $S_{t}$
respectively. Obviously $S_{t}$ is a graph over the minimal 
surface $\Sigma$ if $\Theta>0$ on $S_{t}$.
The evolution equations of $u$ and $\Theta$ can be derived as
follows (cf. \cite{EH1991}),
\begin{equation}\label{eq:evolution of height function}
   \ppl{u}{t}=\biginner{\ppl{F}{t}}{N}=(h-H)\Theta
\end{equation}
and
\begin{equation}\label{eq:evolution of gradient function}
   \ppl{\Theta}{t}=\inner{N}{\nabla{}H}+
        (h-H)\inner{\nablabar_{\nu}N}{\nu}\ .
\end{equation}

\begin{lemma}[{\bf Ecker--Huisken}\ \cite{EH1991}]
The height function $u$ on $S_{t}$ also satisfies
\begin{equation}\label{eq:EH1991-(3.2)}
   \ppl{}{t}\,u=\Delta{}u-\div(\nablabar\ell)+h\Theta\ ,
\end{equation}
where $\div$ is the divergence on $S_{t}$ and $\nablabar$ is the
gradient on $M$.
\end{lemma}

\begin{proof}Since $ u=\ell|{S_{t}}$, we have
$\nabla{}u=(\nablabar\ell)^{\parallel}=\nablabar\ell-\Theta\nu$,
then we obtain
\begin{equation*}
   \Delta{}u=\div\nabla{}u
            =\div(\nablabar\ell)-(\div\nu)\Theta=
            \div(\nablabar\ell)-H\Theta\ .
\end{equation*}
Plugin the above identity to 
\eqref{eq:evolution of height function}, we get 
\eqref{eq:EH1991-(3.2)}.
\end{proof}

\begin{lemma}[{\bf Bartnik}\ \cite{B1984}]
The gradient function $\Theta$ on $S_{t}$ satisfies
\begin{equation}\label{eq:bartnik1984(CMP)-(2.9)}
   \Delta{}\Theta=-(|A|^{2}+\Ricbar(\nu,\nu))\Theta+
   \inner{N}{\nabla{}H}-N(H_{N})\ ,
\end{equation}
where $N(H_{N})$ is the variation of mean curvature of $S_{t}$
under the deformation vector field $N$, which satisfies
\begin{equation}\label{eq:bartnik1984(CMP)-(2.10)}
\begin{aligned}
   N(H_{N})
      =&\,\frac{1}{2}(\nablabar_{\nu}
          \Lcal_{N}\gbar)(e_{i},e_{i})-
          (\nablabar_{e_{i}}\Lcal_{N}\gbar)(\nu,e_{i})-
          \frac{1}{2}\,H\Lcal_{N}\gbar(\nu,\nu)\\
       &\,-\Lcal_{N}\gbar(e_{i},e_{j})\cdot{}A(e_{i},e_{j})\ ,
\end{aligned}
\end{equation}
here $\Lcal$ denotes the Lie derivative.
\end{lemma}

By \eqref{eq:evolution of gradient function} and 
\eqref{eq:bartnik1984(CMP)-(2.9)}, we have the following 
evolution for the gradient function.

\begin{corollary}[{\bf Ecker--Huisken}\ \cite{EH1991}]
$\Theta$ satisfies the following evolution equation
\begin{equation}\label{eq:evolution for Theta}
   \ppl{\Theta}{t}=\Delta{}\Theta+
   (|A|^{2}+\Ricbar(\nu,\nu))\Theta+
   N(H_{N})+(h-H)\inner{\nablabar_{\nu}N}{\nu}\ ,
\end{equation}
where $\Delta$ is the Laplacian on $S_{t}$.
\end{corollary}

Next we will prove that $\{S_{t}\}_{0\leq{}t<T}$ are 
contained in a bounded domain of $M$ for all $T>0$, i.e the
height function is uniformly bounded. This 
result is very important for us to prove 
Theorem \ref{thm:long time existence}. At first, wee need the 
well known maximum principle.

\begin{lemma}[{\bf Maximum Principle}] Let $\Sigma_{1}$ and 
$\Sigma_{2}$ be two hypersurfaces in a Riemannian manifold,
and intersect at a common point tangentially. If $\Sigma_{2}$ 
lies in positive side of $\Sigma_{1}$ around the
common point, then $H_{1}<H_{2}$, where $H_{i}$ is the mean 
curvature of $\Sigma_{i}$ at the common point for $i=1,2$. 
\end{lemma}

\begin{proposition}\label{prop:bound of height function}
Soppose the volume preserving mean curvature
flow \eqref{eq:Vol-MCF} has a family of solutions on $[0,T)$,
$0<T\leq\infty$, then $u$ is uniformly bounded on 
$S\times[0,T)$, i.e.,
\begin{equation*}
   0<C_{1}\leq{}u(x,t)\leq{}C_{2}<\infty\ ,
   \quad\forall\,(x,t)\in{}S\times[0,T)\ ,
\end{equation*}
where $C_{1}$ and $C_{2}$ are two constants depending only
on the initial data $S_{0}(r)=\Sigma(r)$.
\end{proposition}

\begin{proof}At each time $t\in[0,T)$, let $x(t)\in{}S$ be
the point such that
\begin{equation*}
   u_{\max}(t)\equiv\max_{x\in{}S}u(x,t)=u(x(t),t)\ ,
\end{equation*}
and let $y(t)\in{}S$ be the point such that
\begin{equation*}
   u_{\min}(t)\equiv\min_{y\in{}S}u(y,t)=u(y(t),t)\ .
\end{equation*}
Since $\Theta=\inner{N}{\nu}=1$ at $F(x(t),t)$, we have
\begin{equation*}
   0\leq\ppl{u}{t}=h-H\ .
\end{equation*}
By the maximum principle, we have
\begin{equation*}
   h(t)\geq{}H(x(t),t)\geq
   \frac{2\tanh(u_{\max}(t))(1-\Lambda_{+})}
   {1-\tanh^{2}(u_{\max}(t))\Lambda_{+}}\ ,
\end{equation*}
where $\Lambda_{+}=\max\limits_{p'\in{}\Sigma}\lambda^{2}(p')$.
Simlarly, at the point $F(y(t),t)$, we have
\begin{equation*}
   h(t)\leq{}H(y(t),t)\leq
   \frac{2\tanh(u_{\min}(t))(1-\Lambda_{-})}
   {1-\tanh^{2}(u_{\min}(t))\Lambda_{-}}\ ,
\end{equation*}
where $\Lambda_{-}=\min\limits_{p'\in{}\Sigma}\lambda^{2}(p')$.
Thererfore, we have the inequality
\begin{equation*}
   \frac{2\tanh(u_{\min}(t))(1-\Lambda_{-})}
   {1-\tanh^{2}(u_{\min}(t))\Lambda_{-}}
   \geq{}h(t)\geq%
   \frac{2\tanh(u_{\max}(t))(1-\Lambda_{+})}
   {1-\tanh^{2}(u_{\max}(t))\Lambda_{+}}\ .
\end{equation*}
As $t\to{}T$, we have fives cases:
\begin{enumerate}
   \item $u_{\min}(t)\to{}0$ and $u_{\max}(t)\to{}0$;
   \item $u_{\min}(t)\to{}+\infty$ and $u_{\max}(t)\to{}+\infty$;
   \item $u_{\min}(t)\to{}0$ and $u_{\max}(t)\to{}+\infty$;
   \item $u_{\min}(t)$ is uniformly bounded, while
         $u_{\max}(t)\to{}+\infty$;
   \item $u_{\min}(t)\to{}0$, while $u_{\max}(t)$ is
         uniformly bounded.
\end{enumerate}
Case (i) and (ii) could not happen, since the mean curvature flow
is volume preserving. Case (iii) could not happen, otherwise we
would get $0\geq{}2$, a contradiction. Similarly, 
Case (iv) and (v) could not happen.

So the mean curvature flow is uniformly bounded by two surfaces
$\Sigma(r_{1})$ and $\Sigma(r_{2})$ with
$0<r_{1}\leq{}r_{2}<+\infty$ on the time interval $[0,T)$.
\end{proof}

The proof in Proposition \ref{prop:bound of height function}
actually contains the following statement.

\begin{corollary}The average mean curvature $h$ is uniformly
bounded on $[0,T)$, i.e.
\begin{equation*}
   0<\frac{2\tanh(r_{2})(1-\Lambda_{+})}
   {1-\tanh^{2}(r_{2})\Lambda_{+}}
   \leq{}h(t)\leq%
   \frac{2\tanh(r_{1})(1-\Lambda_{-})}
   {1-\tanh^{2}(r_{1})\Lambda_{-}}<2\ .
\end{equation*}
\end{corollary}

\begin{lemma}The mean curvature flow \eqref{eq:Vol-MCF} with
initial data $S_{0}(r)=\Sigma(r)$ preserves the positivity 
of mean curvature of $S_{t}$.
\end{lemma}

\begin{proof}Let
\begin{equation*}
   E(t)=\{x\in{}S\ |\ H(x,t)<0\}
   \quad\text{and}\quad
   E_{t}=F(\cdot,t)(S)\ ,
\end{equation*}
then we have
\begin{equation*}
   \ddl{}{t}\,|E_{t}|=
   -\int_{E_{t}}H(H-h)d\mu<0\ ,
   \quad\forall\,t\in[0,T)\ ,
\end{equation*}
where $|E_{t}|$ denotes the area of $E_{t}$ with respect to
the induced metric $g(t)$ on $S_{t}$,
so $|E_{t}|$ is decreasing. Since $E_{0}=\emptyset$, we
know that $E_{t}=\emptyset$ on $[0,T)$. So the mean curvature of
$S_{t}$ is positive on $[0,T)$.
\end{proof}

Next we will prove that the gradient function $\Theta$ is 
uniformly bounded from below and $|\nabla\Theta|$ is uniformly 
bounded from above on $S_{t}$ for $t\in[0,T)$. 

\begin{proposition}\label{prop:bound for gradient function}
Soppose the volume preserving mean curvature
flow \eqref{eq:Vol-MCF} has a solution on $[0,T)$,
$0<T\leq\infty$, then there exists constants $0<\Theta_{0}<1$
and $0<C_{3}<\infty$ depending only on $S_{0}(r)$ such that
\begin{equation*}
   \Theta\geq{}\Theta_{0}
   \quad\text{and}\quad
   |\nabla\Theta|^{2}\leq{}C_{3}
\end{equation*}
on $S_{t}$ for $0\leq{}t<T$.
\end{proposition}

\begin{proof}Since $\Theta(\cdot,0)\equiv{}1$, we may assume that
$\Theta>0$ for a short time. For any point $p\in{}S_{t}$, 
we may write
\begin{equation*}
   p=(p',u)=(p_{1},p_{2},u)\ ,
\end{equation*}
where $p'=(p_{1},p_{2})\in\Sigma$ and $u$ is the height function on
$S_{t}$. Consider the Gaussian coordinates in $U\times\R\subset{}M$,
where $U\subset\Sigma$ is a neighborhood of $p'$. The unit normal
$\nu$ to $S_{t}$ is given by (cf. \cite[Lemma 3.2]{H1986})
\begin{equation*}
   \nu=\frac{1}{\sqrt{1+|\nabla{}u|^{2}}}
   \left(-\ppl{u}{p_{1}}\,,\,
   -\ppl{u}{p_{2}}\,,\,1\right)\ ,
\end{equation*}
and then the gradient function $\Theta$ is given by
\begin{equation}\label{eq:relation between Theta and u}
   \Theta=\inner{N}{\nu}=\frac{1}{\sqrt{1+|\nabla{}u|^{2}}}\ ,
\end{equation}
where $N=(0,0,1)$. We can see that $|\nabla{}u|=\infty$ if 
and only if $\Theta=0$.

Next, we consider the quasi-linear parabolic equation
\begin{equation}\label{eq:quasilinear equation for u}
   \left\{
   \begin{aligned}
      \ppl{u}{t}&=\Delta{}u-\div(\nablabar\ell)+h\Theta\\
            u(0)&=r\ .
   \end{aligned}
   \right.
\end{equation}
By our hypothesis, \eqref{eq:quasilinear equation for u} has a
solution for $t\in[0,T)$. 
By Proposition \ref{prop:bound of height function}, 
$u$ is uniformly bounded for $t\in[0,T)$. By the standard 
regularity theory of parabolic equation 
(cf. \cite{Lieberman1996} or \cite[Chapter 6]{LSU1967}), 
there exist constants $K_{l}<\infty$ depending only on $l$
and the initial surface $S_{0}(r)$ such that 
\begin{equation*}
      |\nabla^{l}u|\leq{}\ K_{l}, \quad{}l=1,2,\ldots\ , 
\end{equation*}
for $t\in[0,T)$.   

Using \eqref{eq:relation between Theta and u}, these estimates 
imply that $\Theta$ is uniformly
bounded from below and $|\nabla\Theta|^{2}$ is uniformly
from above for $t\in[0,T)$.
\end{proof}

\begin{proposition}\label{prop:bound of A}
Soppose the volume preserving mean curvature
flow \eqref{eq:Vol-MCF} has a family of solutions on $[0,T)$,
$0<T\leq\infty$, then there exists a constant $C_{0}<\infty$ 
depending only on $S_{0}(r)$ such that
\begin{equation*}
   |A|^{2}\leq{}C_{0}<\infty
\end{equation*}
on $S_{t}$ for $0\leq{}t<T$.
\end{proposition}

\begin{proof}We will show that $|A|^{2}$ is uniformly bounded by
contradiction. Let 
$f_{\sigma}=\dfrac{|A|^{2}}{\Theta^{2+\sigma}}$,
where $\sigma>0$ is a small constant. The evolution
equation of $f_{\sigma}$ is given by
\begin{align*}
   \ppl{f_{\sigma}}{t}
      =&\,\Delta{}f_{\sigma}+\frac{2(2+\sigma)}{\Theta}
          \inner{\nabla{}f_{\sigma}}{\nabla\Theta}-
          \frac{2}{\Theta^{2+\sigma}}\,|\nabla{}A|^{2}\\
       &\,+\frac{(1+\sigma)(2+\sigma)|A|^{2}}
           {\Theta^{4+\sigma}}\,|\nabla\Theta|^{2}\\
       &\,+\frac{1}{\Theta^{2+\sigma}}\bigg\{
          -\sigma|A|^{2}(|A|^{2}-2)-2h\tr{}A^{3}
          +8|A|^{2}+2H(h-2H)\\
       &\,\left.-\frac{(2+\sigma)|A|^{2}}{\Theta}\,N(H_{N})
          +\frac{(2+\sigma)|A|^{2}(h-H)}{\Theta}
          \inner{\nablabar_{\nu}N}{\nu}\right\}\ .
\end{align*}

Recall that the restriction to $TS_{t}$ of any tensor field 
$\Phi$ of order $m$ on $M$ can be estimated by
\begin{equation*}
   \|\Phi|_{TS_{t}}(x)\|\leq\Theta^{m}(x)\|\Phi(x)\|\ ,
\end{equation*}
where $\|\Phi(x)\|=|\Phi(x)|_{\gbar_{\alpha\beta}}$ (cf. 
\cite{EH1991}). By using \eqref{eq:bartnik1984(CMP)-(2.10)} 
we estimate the expression $N(H_{N})$ in the evolution equation
\eqref{eq:bartnik1984(CMP)-(2.9)} by
\begin{equation}
   |N(H_{N})|\leq{}C_{4}(\Theta^{3}+\Theta^{2}|A|)\ .
\end{equation}
Here $C_{3}$ depends on $\|\Lcal_{N}\gbar\|_{1}$ where
$\Lcal_{N}\gbar$ is the Lie derivative of the metric with respect to
$N$ whose $C^{1}$-norm can be controlled in terms of $\|N\|_{2}$
(cf. \cite{E2003}). Besides we also have the following estimate
\begin{equation}
   |\inner{\nablabar_{\nu}N}{\nu}|\leq{}C_{5}\Theta^{2}\ ,
\end{equation}
where $C_{5}=\|\nablabar{}N\|$. Since $\{S_{t}\}_{0\leq{}t<T}$ are
contained in a bounded domain whose boundary is 
$\Sigma(r_{1})\cup\Sigma(r_{2})$, the constants $C_{4}$ and
$C_{5}$ only depend on $S_{0}(r)$.

Now assume $|A|_{\max}(t)\to\infty$ as $t\to{}T$. Let 
\begin{equation}
   f_{\max}(t)=\displaystyle\max_{S_{t}}f_{\sigma}\ ,
   \quad\forall\,t\in[0,T)\ .
\end{equation}
Obviously $f_{\max}(t)\geq{}|A|^{2}_{\max}(t)$, 
so $f_{\max}(t)\to\infty$ as $t\to\infty$.
There exists $T_{0}\in(0,T)$ such that when $t>T_{0}$ 
we have the estimate 
\begin{align*}
    \ddl{}{t}\,f_{\max}
        \leq&\,-\sigma\Theta_{0}^{2+\sigma}f_{\max}^{2}
               +(4\sqrt{2}+(2+\sigma)(C_{4}+
                \sqrt{2}\,C_{5}))\Theta_{0}^{1+\sigma/2}
                f_{\max}^{3/2}\\
            &\,+\left(2\sigma+8+(2+\sigma)(C_{4}+2C_{5})+
                \frac{(1+\sigma)(2+\sigma)C_{3}}
                {\Theta_{0}^{2}}\right)f_{\max}\\
        \leq&\,-\frac{\sigma\Theta_{0}^{2+\sigma}}{2}
                f_{\max}^{2}\ .
\end{align*}
This is a contradiction since $df_{\max}/dt\geq{}0$. 
Therefore $f_{\sigma}$ must be 
uniformly bounded, which implies that $|A|^{2}$ must be 
uniformly bounded.
\end{proof}

\begin{proposition}[{{\bf Huisken}\ \cite[$\S${}4]{H1987}}]
\label{prop:H1987-(4.1)}
For every natural number $m$, we have the following 
evolution equation:
\begin{equation}
   \begin{aligned}
      \ppl{}{t}\,|\nabla^{m}A|^{2}
          =&\,\Delta|\nabla^{m}A|^{2}-2|\nabla^{m+1}A|^{2}+
              \sum_{i+j+k=m}\nabla^{i}A*\nabla^{j}A*
              \nabla^{k}A*\nabla^{m}A\\
           &\,+h\sum_{i+j=m}\nabla^{i}A*\nabla^{j}A*\nabla^{m}A\ .
   \end{aligned}
\end{equation}
Furthermore, there exists constant $\{C(m)\}_{m=1,2,\ldots}$ 
depending only on $m$ and $S_{0}(r)$ such that
\begin{equation}
   |\nabla^{m}A|^{2}\leq{}C(m)\ ,
   \quad{}m=1,2,\ldots\ ,
\end{equation}
are uniformly on $S_{t}$ for $0\leq{}t<T$. 
\end{proposition}

By the above discussion, the constants in 
Proposition \ref{prop:bound of height function} and Proposition
\ref{prop:bound for gradient function}--
\ref{prop:H1987-(4.1)} are independent of time. 
Now we can prove part one of Theorem~\ref{thm:long time existence}.

\begin{proof}[{\bf Proof of Theorem~\ref{thm:long time existence}}]
(1) (cf. \cite{H1984,H1987}) Assume that $T<\infty$. Let
\begin{equation}\label{eq:limit surface ST}
   S_{T}=\lim_{t\to{}T}S_{t}
        =\left\{\lim_{t\to{}T}F(x,t)\ \bigg|\ x\in{}S\right\}\ .
\end{equation}
We claim that $S_{T}$ is a smooth surface which is 
homeomorphic to $S$.

In fact, by Proposition \ref{prop:bound of height function}, 
the height function $u$ is uniformly bounded on $S_{t}$ for 
$t\in[0,T)$. So \eqref{eq:limit surface ST} is well defined. 
Since $|A|^{2}$ is uniformly bounded for $t\in[0,T)$, we have
\begin{equation*}
   \int_{0}^{T}\max_{S_{t}}\left|\ppl{}{t}\,g_{ij}\right|dt
   \leq{}C<\infty\ ,
\end{equation*}
so $S_{T}$ is a well defined surface by Lemma 14.2 in 
\cite{Hamilton1982}. Since $|\nabla^{m}A|^{2}$, $m=1,2,\ldots$, 
are uniformly bounded for $t\in[0,T)$, $S_{T}$ is smooth. 

Now we consider a new volume preserving mean curvature flow
\begin{equation*}
   \ppl{F}{t}=(h-H)\nu
\end{equation*}
with initial data $S_{T}$. This flow has a short time solution
for $t\in[T,T_{1})$, where $T_{1}>T$, the detail can be found in
\cite[$\S${}6.7]{CK2004}. This contradicts to the assumption 
that $T$ is maximal. Therefore the maximal time $T$ of the
volume preserving mean curvature flow \eqref{eq:Vol-MCF} must be
infinite.
\end{proof}

\subsection{Exponential convergence to CMC surfaces} 
We have proved that the volume preserving mean curvature flow 
\eqref{eq:Vol-MCF} has a long time solution. Let 
\begin{equation}\label{eq:limiting surface}
   S_{\infty}(r)=\lim_{t\to\infty}S_{t}
\end{equation} 
be the limiting surface. Obviously $S_{\infty}(r)$ has 
the following properties: 
\begin{enumerate}
   \item It is well defined since $\{S_{t}\}_{0\leq{}t<\infty}$ 
         are contained in a bounded domain of $M$. 
   \item It's also a smooth surface since $|\nabla^{m}A|^{2}$,
         $m=0,1,2,\ldots$, are uniformly bounded for 
         $t\in[0,\infty)$.
   \item It's a graph over $\Sigma$ since $\Theta$ is uniformly
         bounded from below for $t\in[0,\infty)$.
\end{enumerate}
In this subsection, we will show that the solution 
surface $S_{t}$ converges exponentially fast to $S_{\infty}(r)$ 
(cf. \cite{CrM2007,H1987,HY1996}), although we don't need this
fact to prove the existence of the CMC foliation of $M$.

\begin{proposition}Suppose $(S_{t},g(t))$ is a solution to the
mean curvature flow \eqref{eq:Vol-MCF} for $t\in[0,\infty)$, 
then
\begin{equation}
   \lim_{t\to\infty}\sup_{S_{t}}|H-h|=0\ .
\end{equation}
Therefore $S_{\infty}(r)$ is a surface of constant mean curvature.
\end{proposition}

\begin{proof}Since
\begin{equation*}
   \ddl{}{t}\,|S_{t}|=-\int_{S_{t}}(H-h)^{2}d\mu<0\ ,
\end{equation*}
where $|S_{t}|$ denotes the area of $S_{t}$ with respect to the 
metric $g(t)$, then we have
\begin{equation*}
   \int_{0}^{\infty}\int_{S_{t}}(H-h)^{2}d\mu{}dt
   \leq{}|S_{0}|\ .
\end{equation*}
On the other hand, by Lemma~\ref{lem:HY1996-(3.6)} and 
Lemma~\ref{lem:HY1996-(3.7)}, we have
\begin{align*}
   \ddl{}{t}\int_{S_{t}}(H-h)^{2}d\mu
         =&\,2\int_{S_{t}}(H-h)\ddl{}{t}(H-h)d\mu-
              \int_{S_{t}}H(H-h)^{3}d\mu\\
         =&\,2\int_{S_{t}}(H-h)[\Delta{}H+(H-h)(|A|^{2}-2)]d\mu\\
          &\,-\int_{S_{t}}H(H-h)^{3}d\mu\\
         =&\,-2\int_{S_{t}}|\nabla{}H|^{2}d\mu+
              2\int_{S_{t}}(H-h)^{2}(|A|^{2}-2)d\mu\\
          &\,-\int_{S_{t}}H(H-h)^{3}d\mu\ ,
\end{align*}
here we use the identity $\int_{S_{t}}(H-h)d\mu=0$.
By Proposition~\ref{prop:H1987-(4.1)} and the inequalities
$|\nabla{}H|\leq\sqrt{2}\,|\nabla{}A|$, there is a constant
$C_{6}<\infty$ depending only on $S_{0}(r)$ such that 
\begin{equation}\label{eq:2th derivative of |S_{t}|}
   \left|\ddl{}{t}\int_{S_{t}}(H-h)^{2}d\mu\right|
   \leq{}C_{6}
\end{equation}
is uniformly for $t\in[0,\infty)$. So we have
\begin{equation}
   \lim_{t\to\infty}\int_{S_{t}}(H-h)^{2}d\mu=0\ .
\end{equation}
Then for any $p>2$, by the interpolation arguments 
(cf. \cite[$\S${}5]{CrM2007} for detail), the inequality
$|\nabla^{2}H|\leq\sqrt{2}\,|\nabla^{2}A|$ and 
Proposition~\ref{prop:H1987-(4.1)}, we have
\begin{align*}
   \sup_{S_{t}}|H-h|
       &\leq{}C\|\nabla^{2}H\|_{2}^{1/p}\|H-h\|_{2}^{1/p}\\
       &\leq{}C\left(\int_{S_{t}}(H-h)^{2}
        d\mu\right)^{1/(2p)}\\
       &\to{}0\quad(\text{as}\ t\to\infty)\ .
\end{align*}
where $\|\cdot\|_{2}=\|\cdot\|_{L^{2}(S_{t})}$.  
So the proposition follows.
\end{proof}

We say that a surface $S$ with constant mean curvature is 
(strictly) {\em stable} if volume preserving variations of $S$ 
in $M$ incease the area, or equivalently if the second variation
operator on $S$,
\begin{equation*}
   L\phi=-\Delta\phi-(|A|^{2}+\Ricbar(\nu,\nu))\phi
\end{equation*}
has only strictly positive eigenvalues when restricted to 
functions $\phi$ with $\displaystyle\int_{S}\phi\,d\mu=0$.

\begin{lemma}\label{lem:stability of limit surface}
For each $r\in\R$, the limit surface $S_{\infty}(r)$
to the volume preserving mean curvature flow \eqref{eq:Vol-MCF}
is strictly stable surface of constant mean curvature. 
\end{lemma}

\begin{proof}Suppose $S'$ is a volume preserving variation of 
$S_{\infty}(r)$, such that $S'$ is a graph over $\Sigma$ and
$\Area(S')<\Area(S_{\infty}(r))$. Consider the volume preserving
mean curvature flow \eqref{eq:Vol-MCF} with initial surface $S'$.
By the above discussion, there is a long time solution to this
volume preserving mean curvature flow. 
Let $S_{\infty}'$ be the limiting surface, then it is a graph over
$\Sigma$ whose mean curvature is a constant and 
$\Area(S_{\infty}')<\Area(S_{\infty}(r))$. 

We claim that this is impossible. In fact, according to Theorem 
\ref{thm:main theorem}, $\{S_{\infty}(r)\}_{r\in\R}$ foliate $M$, 
so there are two surfaces $S_{\infty}(r_{1})$ and 
$S_{\infty}(r_{2})$, where $r_{1}<r_{2}$, which
touch $S'$ from the below and from the above for the first time
respectively. By maximum principle, we have
\begin{equation*}
   H(S_{\infty}(r_{2}))<H(S_{\infty}')<H(S_{\infty}(r_{1}))\ .
\end{equation*}
But this is impossible since 
$H(S_{\infty}(r_{1}))<H(S_{\infty}(r_{2}))$ when $r_{1}<r_{2}$
(see the proof of Theorem \ref{thm:main theorem} in 
$\S${}\ref{sec:foliation}). So the stability of limiting surfaces
follows.
\end{proof}

\begin{proof}[{\bf Proof of Theorem~\ref{thm:long time existence}}]
(2) Since $S_{\infty}(r)$ is stable, the lowest eigenvalue 
$\lambda_{\infty}$ of the Jacobi operator $L_{\infty}$ on 
$S_{\infty}(r)$ is positve, where 
\begin{equation*}
   L_{\infty}\phi=-\Delta_{\infty}\phi-
   (|A_{\infty}|^{2}-2)\phi\ ,
\end{equation*}
here $\Delta_{\infty}$ is the Laplacian on $S_{\infty}(r)$ and 
$A_{\infty}$ is the second fundamental form of $S_{\infty}(r)$.
Let $\lambda_{t}$ be the lowest eigenvalue of the Jacobi operator
$L$ on $S_{t}$. Then $\lambda_{t}\to\lambda_{\infty}$ 
as $t\to\infty$. For any 
$0<\varepsilon<\dfrac{2}{3}\,\lambda_{\infty}$, 
there exists $T>0$ such that for any $t>T$ we have
\begin{equation*}
   |\lambda_{\infty}-\lambda_{t}|<\varepsilon 
   \qquad\text{and}\qquad
   \sup\limits_{S_{t}}|H(H-h)|\leq\varepsilon\ .
\end{equation*}
Therefore, when $t>T$ we have
\begin{align*}
   \ddl{}{t}\int_{S_{t}}(H-h)^{2}d\mu\leq{}
   -(2\lambda_{\infty}-3\varepsilon)\int_{S_{t}}(H-h)^{2}d\mu\ ,
\end{align*}
which implies 
\begin{equation*}
    \int_{S_{t}}(H-h)^{2}d\mu\leq
    \left(\int_{S_{T}}(H-h)^{2}d\mu\right)
    e^{-(2\lambda_{\infty}-3\varepsilon)t}\ .
\end{equation*}
By the same interpolation arguments as above, we know that
$\sup|H-h|$ converges exponentially to zero. Since
\begin{equation*}
   \left|\ppl{F}{t}\right|=|h-H|\ ,
\end{equation*}
we obtain that $S_{t}$ converges exponentially to the limiting
surface which has constant mean curvature.
So Part two of Theorem \ref{thm:long time existence} is proved.
\end{proof}

\section{Existence of CMC foliation}\label{sec:foliation}

We need a lemma of Mazzeo and Pacard  which will be useful for 
proving the uniqueness of the CMC foliation of $M$.

\begin{lemma}[{\bf Mazzeo--Pacard}\ \cite{MP2007}]
\label{lem:MP2007-(3.1)}
Suppose that $\Fcal$ is a monotonically increasing CMC foliation in
$(M,\gbar_{\alpha\beta})$, then $\Fcal$ is unique amongst all CMC
foliations whose leaves are diffeomorphic to $\Sigma$.
\end{lemma}

\begin{proof}[{\bf Proof of Theorem \ref{thm:main theorem}}]
(1) At first, we can foliate the quasi-Fuchsian $3$-manifold $M$ by 
the surfaces $\Sigma(r)$, $r\in\R$. All of these surfaces, except
$\Sigma\equiv\Sigma(0)$ (the minimal surface), are not surfaces of
constant mean curvature. But for each $r>0$, we consider the
mean curvature flow \eqref{eq:Vol-MCF} with initial condition
$S_{0}=\Sigma(r)$. By Theorem \ref{thm:long time existence}, 
we have a solution of \eqref{eq:Vol-MCF}, which is a smooth 
surface of (positive) constant mean curvature, and we denote it 
by $S_{\infty}(r)$. For these surfaces $\Sigma(r)$ with $r<0$, 
we have the surfaces with (negative) constant mean curvature.
We need three steps to prove that the limiting surfaces
$S_{\infty}(r)$, $r\in\R$, form a CMC foliation of $M$. 

{\bf Step 1}: {\em The limiting surfaces are embedded.} 
This is obviously since each surface $S_{\infty}(r)$ is a graph 
over the minimal surface $\Sigma$.

{\bf Step 2}: {\em The limiting surfaces are disjoint.}  
Assume that $0<r_{1}<r_{2}$, we will show that 
$S_{\infty}(r_{1})\cap{}S_{\infty}(r_{2})=\emptyset$. Consider
two volume preserving mean curvature flows \eqref{eq:Vol-MCF}
with initial data $\Sigma(r_{1})$ and $\Sigma(r_{2})$ respectively.
Let $u_{1}$ and $u_{2}$ be the height functions of the surfaces
$S_{t}(r_{1})$ and $S_{t}(r_{2})$ respectively, then we have
$u_{1}(x,0)<u_{2}(x,0)$ for all $x\in{}S$. Now we assume 
that two surfaces $S_{t}(r_{1})$ and $S_{t}(r_{2})$ touch 
for the first time at $T_{0}\in(0,\infty)$ and $p_{0}\in{}M$. 
Recall that the height functions satisfy the evolution equation
\eqref{eq:EH1991-(3.2)}. Let $w=u_{2}-u_{1}$, then $w\geq{}0$, 
and around $p_{0}$ we have
\begin{equation*}
   0>Lw=\Delta{}w+\inner{\cdot}{\nabla{}w}-\ppl{w}{t}\ ,
\end{equation*}
here we use the fact that $h_{1}(t)<h_{2}(t)$ since 
$H(S_{t}(r_{1}))<H(S_{t}(r_{2}))$ pointwise, where $h_{1}(t)$ 
and $h_{2}(t)$ are the average mean curvature of $S_{t}(r_{1})$ 
and $S_{t}(r_{2})$ respectively. By the strong maximum 
principle (cf. \cite{Friedman1964,PW1967}), this is impossible 
unless $w\equiv{}0$. But $w\equiv{}0$ implies $u_{1}\equiv{}u_{2}$,
which is also impossible since the flows preserve volume. 
This means that $S_{t}(r_{1})$ and $S_{t}(r_{2})$ are
disjoint all the time, so $S_{\infty}(r_{1})$ and
$S_{\infty}(r_{2})$ are disjoint.

{\bf Step 3}: {\em We claim}
\begin{equation*}
   M=\bigcup_{r\in\R}S_{\infty}(r)\ .
\end{equation*} 
In fact, according to the proof of 
Proposition \ref{prop:bound for gradient function}, for each
$r\ne{}0$, $\Sigma\cap{}S_{\infty}(r)=\emptyset$. Let
$Q(r)$ be the domain bounded by $\Sigma$ and $S_{\infty}(r)$.
Since $\{\Sigma(r)\}_{r\in\R}$ foliate $M$ and each
$S_{\infty}(r)$ is the limiting surface of the volume 
preserving mean curvature flow with initial data $\Sigma(r)$,
the volume of $Q(r)$ is a continuous function with respect to $r$. 
Together with the facts that the limiting surfaces are 
embedded and disjoint, Step 3 is proved.

Therefore these surfaces form a CMC foliation of $M$.

(2) We claim that the foliation 
$\Fcal=\{S_{\infty}(r)\}_{r\in\R}$
is monotonically increasing: if $r_{1}<r_{2}$, then
$H(S_{\infty}(r_{1}))<H(S_{\infty}(r_{2}))$. In fact, 
since $H$ satisfies the (strictly) parabolic equation: 
\begin{equation*}
   \ppl{H}{t}=\Delta{}H+(H-h)(|A|^{2}-2)\ ,
\end{equation*}
and $H(\Sigma(r_{1}))<H(\Sigma(r_{2}))$ pointwise, then by 
the comparison principle for quasilinear parabolic equations
(cf.\cite[Theorem 9.7]{Lieberman1996}), we have
$H(S_{t}(r_{1}))<H(S_{t}(r_{2}))$ pointwise for 
$t\in[0,\infty)$. In particular, 
$H(S_{\infty}(r_{1}))<H(S_{\infty}(r_{2}))$.

Since this foliation is monotonically increasing,
we get the uniqueness of the CMC foliation by 
Lemma \ref{lem:MP2007-(3.1)}.
\end{proof}

\begin{remark}In \cite{Toda1999}, M. Toda proved so called 
volume constraint Plateau problem in hyperbolic $3$-manifolds
satisfying some conditions. Our quasi-Fuchisan manifolds satisfy 
the conditions required in his paper, so for each $r\in\R$, we can 
find an area minimizing surface $S(r)$ such that the volume of 
the domain bounded by $\Sigma$ and $S(r)$ is equal to the volume 
of the domain bounded by $\Sigma$ and $\Sigma(r)$. Each $S(r)$ is 
a surface of constant mean curvature. If one can show that
$S(r_{1})\cap{}S_(r_{2})=\emptyset$ for $r_{1}\ne{}r_{2}$ and
$M=\cup{}S_{r}$, then $\{S_{r}\}_{r\in\R}$ form a CMC foliation
of $M$.
\end{remark}

\section{A counterexample}\label{sec:counterexam}

In this section, we will show that Theorem \ref{thm:main theorem}
is not true for the quasi-Fuchsian $3$-manifolds containing
minimal surfaces with big principle curvature.

\subsection{Existence of the surfaces with CMC}
We need some results of J. Gomes and R. L{\'o}pez (cf.
\cite{Gomes1987,Lopez2000}). Let $\H^{3}$ be a 
three-dimensional hyperbolic space of constant
sectional curvature $-1$. We will work in the Poicar{\'e} model
of $\H^{3}$, i.e.,
\begin{equation*}
   \H^{3}=\{(x,y,z)\in\R^{3}\ |\ 
   x^{2}+y^{2}+z^{2}<{}1\}
\end{equation*}
equipped with metric
\begin{equation*}
   ds^{2}=\frac{4(dx^{2}+dy^{2}+dz^{2})}{(1-r^{2})^{2}}\ ,
\end{equation*}
where $r=\sqrt{x^{2}+y^{2}+z^{2}}$.
The hyperbolic space $\H^{3}$ has a natural compactification
$\overline{\H}{}^{3}=\H^{3}\cup{}S_{\infty}^{2}$, where
$S_{\infty}^{2}=\widehat{\C}$ is the Riemann sphere. Suppose
$X$ is a subset of $\H^{3}$, we call the set 
$\partial_{\infty}X$ defined by
\begin{equation*}
   \partial_{\infty}X=\overline{X}
   \cap{}S_{\infty}^{2}\ ,
\end{equation*}
the {\em asymptotic boundary} of $X$, where 
$\overline{X}$ is the closure of $X$ in 
$\overline{\H}{}^{3}$.

Suppose $G$ is a subgroup $\Isom(\H^{3})$ which leaves a 
geodesic $\gamma\subset\H^{3}$ pointwise fixed. We call $G$
the spherical group of $\H^{3}$ and $\gamma$ the rotation axis
of $G$. A surface in $\H^{3}$ invariant by $G$ is called a 
{\em spherical surface}. For two circles $C_{1}$ and $C_{2}$ in 
$\H^{3}$, if there is a geodesic $\gamma$ such that each of 
$C_{1}$ and $C_{2}$ is invariant by the group of rotations that 
fixes $\gamma$ pointwise, then $C_{1}$ and $C_{2}$ are said to 
be {\em coaxial}, and $\gamma$ is called the 
{\em rotation axis} of $C_{1}$ and $C_{2}$.

Let $P_{1}$ and $P_{2}$ be two disjoint geodesic plane in 
$\H^{3}$. Then $P_{1}\cup{}P_{2}$ divides $\H^{3}$ in three
components. Let $X_{1}$ and $X_{2}$ be the two of 
them with $\partial{}X_{i}=P_{i}$ for $i=1,2$. Given 
two subsets $A_{1}$ and $A_{2}$ of $\overline{\H}{}^{3}$, we 
say $P_{1}$ and  $P_{2}$ {\em separate} $A_{1}$ and $A_{2}$ if
one of the following cases occurs (cf. \cite{Lopez2000}):
\begin{enumerate}
   \item if $A_{1},A_{2}\subset\H^{3}$, then
         $A_{i}\subset{}X_{i}$ for $i=1,2$;
   \item if $A_{1}\subset\H^{3}$ and 
         $A_{2}\subset{}S_{\infty}^{2}$, then 
         $A_{1}\subset{}X_{1}$ and 
         $A_{2}\subset\partial_{\infty}X_{2}$;
   \item if $A_{1},A_{2}\subset{}S_{\infty}^{2}$, then
         $A_{i}\subset\partial_{\infty}X_{i}$ for $i=1,2$.
\end{enumerate}
Then we may define the distance between $A_{1}$ and $A_{2}$ by
\begin{equation}\label{eq:distance between circles}
   d(A_{1},A_{2})=\sup\{\dist(P_{1},P_{2})\ |\ 
   P_{1}\ \text{and}\ P_{2}\ \text{separate}\ A_{1}\
   \text{and}\ A_{2}\}\ ,
\end{equation}
where $\dist(P_{1},P_{2})$ is the hyperbolic distance between 
$P_{1}$ and $P_{2}$. 

\begin{lemma}[{\bf Gomes}\ \cite{Gomes1987}]
\label{lem:Gomes1987}
There exists a finite
constant $d_{0}>0$ such that for two disjoint circles 
$C_{1},C_{2}\subset{}S_{\infty}^{2}$, if 
$d(C_{1},C_{2})\leq{}d_{0}$, then there exists a minimal surface
$\Pi$ which is a surface of revolution and whose asymptotic
boundary is $C_{1}\cup{}C_{2}$.
\end{lemma}

Let $C_{1}$ and $C_{2}$ be two disjoint circles on
$S_{\infty}^{2}$, and let $P_{1}$ and $P_{2}$ be two geodesic 
planes whose asymptotic boundaries are $C_{1}$ and $C_{2}$
respectively. Suppose $C'_{1}\subset{}P_{1}$ and 
$C'_{2}\subset{}P_{2}$ so that $C'_{1}$ and $C'_{2}$
are two coaxial circles with respect to the rotation axis of
$C_{1}$ and $C_{2}$.

\begin{lemma}[{\bf L{\'o}pez}\ \cite{Lopez2000}]
\label{lem:Lopez2000} 
Given $H\in(-1,1)$, there exists a constant $d_{H}$ depending 
only on $H$ such that if $d(C_{1},C_{2})\leq{}d_{H}$, then 
there exists a surface $\Pi$ contained in the domain bounded 
by $P_{1}$ and $P_{2}$ such that
\begin{itemize}
   \item $\Pi$ is a surface of revolution whose boundary is
         $C'_{1}\cup{}C'_{2}$, and
   \item $\Pi$ is a surface whose mean curvature is equal to
         $H$ with respect to the normal pointing to the domain
         containing the rotation axis of $C_{1}$ and $C_{2}$.
\end{itemize}
\end{lemma}

\begin{remark}In Lemma~\ref{lem:Lopez2000}, when $H<0$, then 
there is no such a surface $\Pi$ if we replace $C'_{i}$ by
$C_{i}$ for $i=1,2$ (cf. \cite{P1999}).
\end{remark}

\subsection{Detail description of the counterexample} 
Now we choose four circles $\{C_{i}\}_{i=1,\ldots,4}$ on 
$S_{\infty}^{2}$ such that
$d(C_{1},C_{2})$ and $d(C_{3},C_{4})$ are sufficiently small, 
where $d(\cdot,\cdot)$ is the distance defined by
\eqref{eq:distance between circles}. Let $D_{i}$ be the geodesic
plane in $\H^{3}$ such that $\partial_{\infty}D_{i}=C_{i}$ for
$i=1,\ldots,4$. By some M\"obius 
transformation, we may assume that the middle point of the
geodesic segment which is perpendicular to both $D_{1}$ and 
$D_{2}$ passes through the origin.

For any circle $C\subset{}S_{\infty}^{2}$, we may define the
distance between the origin $O$ (or any fixed point) and the 
circle $C$ to be the hyperbolic distance between $O$ and the 
geodesic plane whose asymptotic boundary is $C$. Because of this
definition, we may say that the radius of the circle $C$ is 
big or small if the distance between $O$ and $C$ is small or big.
  
Let $\Lambda$ be a closed smooth curve on $S_{\infty}^{2}$, 
then cover $\Lambda$ by finite disks 
$\{B_{l}\subset{}S_{\infty}^{2}\}_{l=1,\ldots,N}$ 
with small radii such that
\begin{itemize}
   \item each circle $\partial{}B_{l}$ is invariant under the 
         rotation along the geodesic connecting the origin $O$ and
         the center of the disk $B_{l}$, which locates at 
         $\Lambda$,
   \item the radii of disks are small enough so that
         $B_{l}\cap{}C_{i}=\emptyset$ for $l=1,\ldots,N$ and 
         $i=1,\ldots,4$, and
   \item for each $l\equiv{}1\ (\mod{}N)$, $\partial{}B_{l}$ 
         intersects both
         $\partial{}B_{l-1}$ and $\partial{}B_{l+1}$ and no
         other circle,
\end{itemize} 
then we get a quasi-Fuchsian group $\Gamma$ which is the subgroup 
of orientation preserving transformations in the group generated 
by $N$ reflections about the circles 
$\partial{}B_{1},\ldots,\partial{}B_{N}$ 
(cf. \cite [Page 263]{Bers1972} or \cite[Page 149]{Bers1981}). 
The limit set of the quasi-Fuchsian group $\Gamma$, denoted 
by $\Lambda_{\Gamma}$, is around the curve $\Lambda$. Let 
$S_{\infty}^{2}\setminus\Lambda_{\Gamma}=\Omega_{1}\cup\Omega_{2}$, 
where $\Omega_{1}$ contains $C_{1}$ and $C_{2}$, while 
$\Omega_{2}$ contains $C_{3}$ and $C_{4}$. See Figure 1.

\begin{figure}[htbp]
\begin{center}
      \includegraphics[scale=0.5]{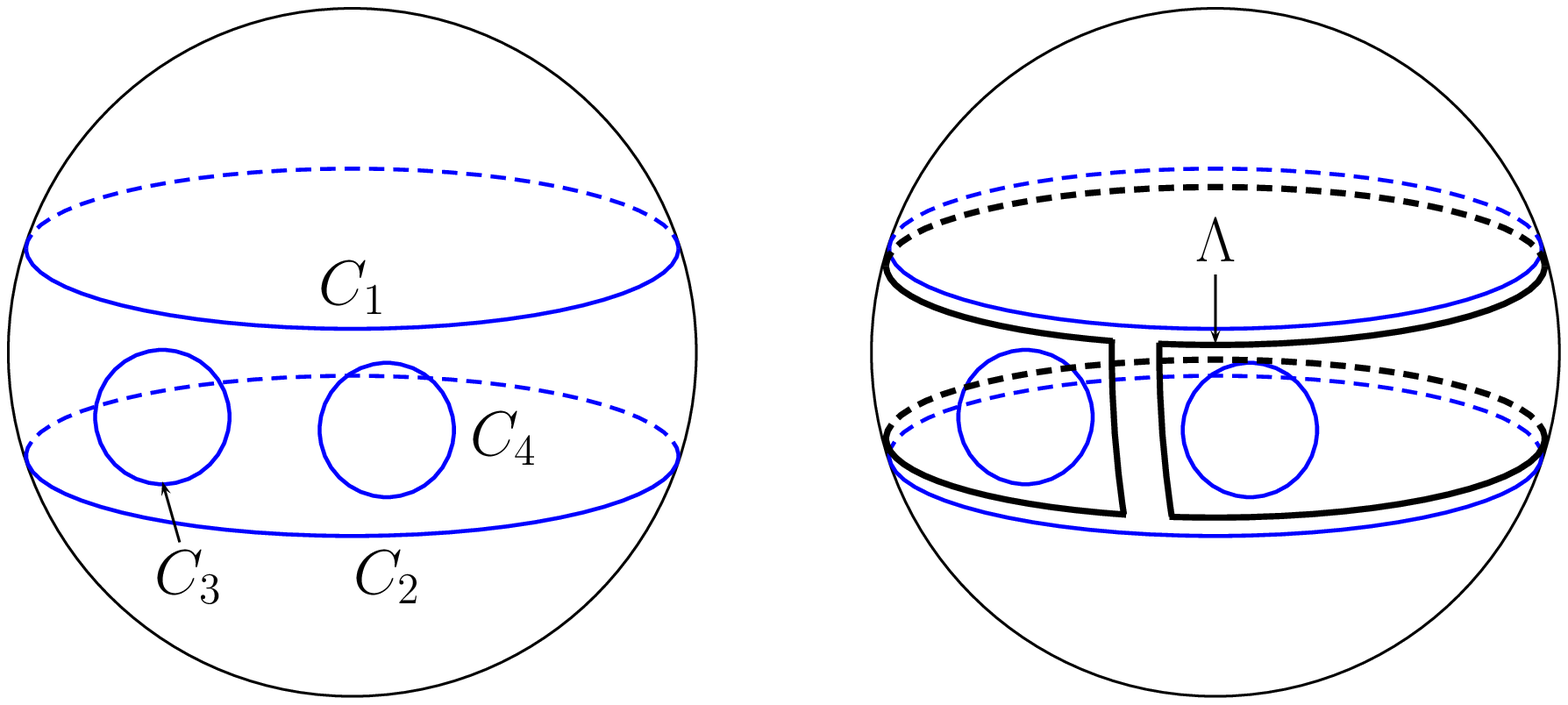}
      \caption{}
\end{center}
\end{figure}

\textbf{Claim}: {\em The quasi-Fuchsian $3$-manifold 
$\H^{3}/\Gamma$ constructed above can not be foliated by 
surfaces of constant mean curvature}.

Let $\varepsilon>0$ be sufficiently small, and 
let $H_{0}=2\tanh\varepsilon$. Let $d_{0}$ and 
$d_{H_{0}}$ be two constants given in Lemma~\ref{lem:Gomes1987} 
and Lemma~\ref{lem:Lopez2000}, and suppose
$d(C_{1},C_{2})=2\varepsilon\ll{}d_{0}$ and 
$d(C_{3},C_{4})\ll{}\min\{d_{H_{0}},d_{0}\}$. 

Now assume that $\H^{3}/\Gamma$ is foliated by surfaces of constant 
mean curvature, where each surface is closed and is homotopic
to $\H^{3}/\Gamma$. Lift the foliation to the universal covering 
space $\H^{3}$, then there should exist a foliation of $\H^{3}$ 
so that each leaf is a disk with constant mean curvature and 
with the same asymptotic boundary $\Lambda_{\Gamma}$.
Notice that any disk type surface in $\H^{3}$ with asymptotic 
boundary $\Lambda_{\Gamma}$ divides $\overline{\H}{}^{3}$ into 
two parts, one of them contains $C_{1}$ and $C_{2}$, while 
the other contains $C_{3}$ and $C_{4}$. 
We choose a normal vector field on the disk type surface so that 
each normal vector points to the domain containing $C_{1}$ and 
$C_{2}$. Assume that there is a CMC foliation $\Fcal=\{L_{t}\}$ 
with a parameter $t\in(-\infty,\infty)$ such that
\begin{itemize}
   \item the leaves are convergent to $\Omega_{1}$ as 
         $t\to{}-\infty$ and
   \item the leaves are convergent to $\Omega_{2}$ as 
         $t\to{}\infty$.
\end{itemize} 
In other words, we have
\begin{equation}\label{eq:foliation limit behavior}
   \lim_{t\to\pm\infty}H(L_{t})=\pm{}2\ ,
\end{equation}
where $H(L_{t})$ denotes the mean curvature of the leaf $L_{t}$
with respect to the normal vector pointing to the domain 
containing $C_{1}$ and $C_{2}$.

Since $d(C_{3},C_{4})$ is very small, there exists a minimal 
surface with asymptotic boundary $C_{3}\cup{}C_{4}$ by 
Lemma~\ref{lem:Gomes1987}. Consider the leaf $L_{t'}\in\Fcal$ which 
touches the minimal surface for the first time, then the mean
curvature of $L_{t'}$ must be positive by the maximal principle. 
Because of \eqref{eq:foliation limit behavior}, there exists 
$-\infty<t_{1}<t'$ such that the mean curvature of $L_{t_{1}}$
is zero, i.e. the leaf $L_{t_{1}}$ is a disk type minimal surface.
Similarly, we have another leaf $L_{t_{2}}\in\Fcal$ which is a 
disk type minimal surface with asymptotic boundary $\Lambda$. See
Figure 2.

\begin{figure}[htbp]
\begin{center}
      \includegraphics[scale=0.5]{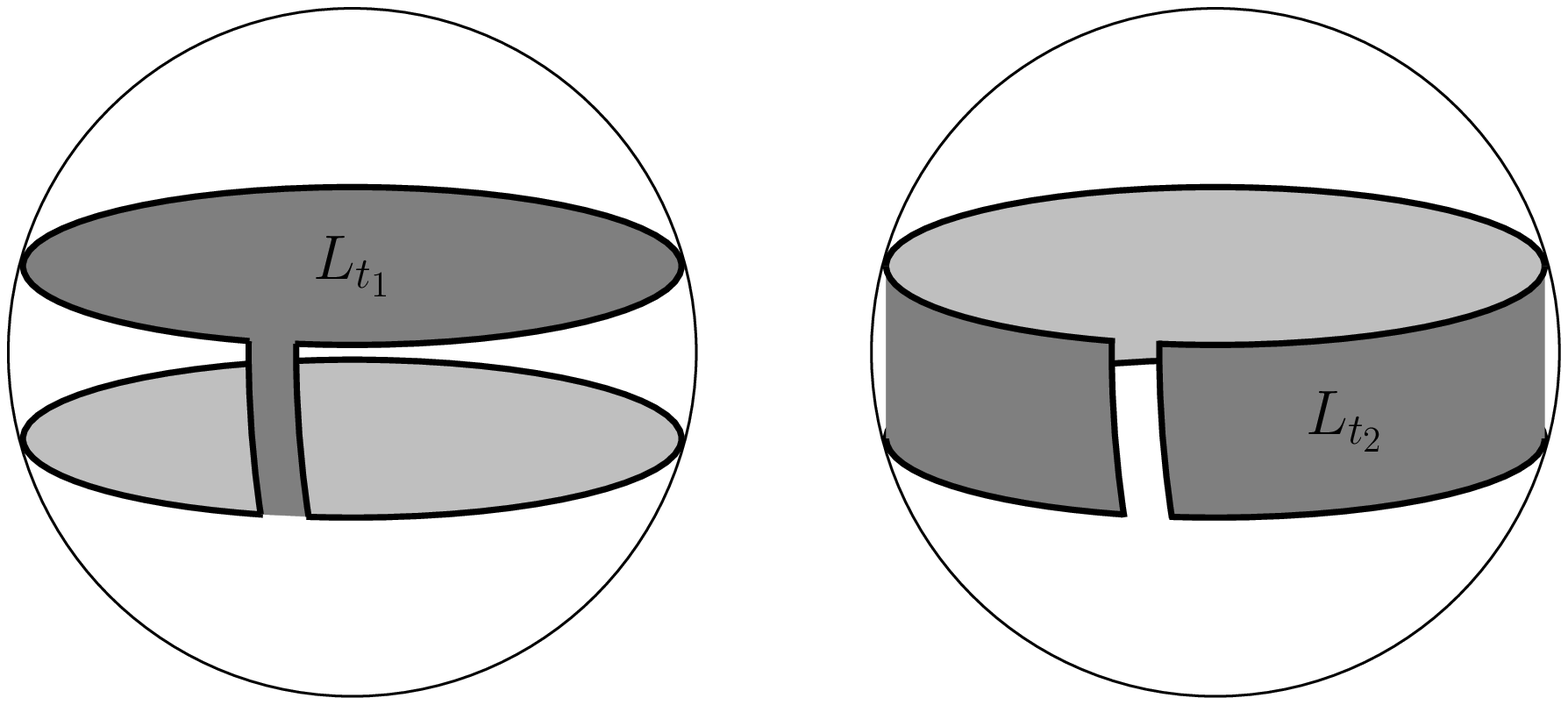}
      \caption{}
\end{center}
\end{figure}

Let $X\subset\H^{3}$ be the domain bounded by $L_{t_{1}}$ and
$L_{t_{2}}$, then by assumption $X$ is foliated by 
$\{L_{t}\}_{t_{1}\leq{}t\leq{}t_{2}}$, i.e.
\begin{equation*}
   X=\bigcup_{t_{1}\leq{}t\leq{}t_{2}}L_{t}\ .
\end{equation*}
Notice that $D_{3}$ and $D_{4}$ are disjoint from $X$. We choose
two circles $C_{3}'\subset{}D_{3}$ and $C_{4}'\subset{}D_{4}$
so that $C_{3}'$ and $C_{4}'$ are coaxial with respect to the
rotation axis of $C_{3}$ and $C_{4}$, by 
Lemma~\ref{lem:Lopez2000} there is a surface $\Pi_{0}$
with constant mean curvature $-H_{0}$ with 
respect to the normal pointing to the domain containing
the rotation axis of $C_{3}'$ and $C_{4}'$. Obviously 
$\Pi_{0}$ is disjoint form $L_{t_{1}}$ but intersects 
$L_{t_{2}}$. Let $\Pi_{0}'=\Pi_{0}\cap{}X$. Consider the leaf 
\begin{equation*}
    L_{t''}\in\{L_{t}\ |\ t_{1}\leq{}t\leq{}t_{2}\}
\end{equation*}
which  touches $\Pi_{0}'$ for the first time, then 
$H(L_{t''})>H_{0}$ by the maximal principle. So 
there exists $t_{3}\in(t_{1},t_{2})$ such that 
$H(L_{t_{3}})=H_{0}$. We claim that the leaf
$L_{t_{3}}$ must self-intersects.

Let $D_{1}(\varepsilon)$ be the disk bounded by $C_{1}$ with
$H(D_{1}(\varepsilon))=H_{0}$ with respect to the
normal vector pointing to domain not containing $C_{2}$, and 
similarly let $D_{2}(\varepsilon)$ be the disk bounded by 
$C_{2}$ with $H(D_{2}(\varepsilon))=H_{0}$ with 
respect to the normal vector pointing to domain not 
containing $C_{1}$. Then 
$D_{1}(\varepsilon)\cap{}D_{2}(\varepsilon)=\{O\}$, where
$O\in\H^{3}$ is the origin. By maximal principle, both 
$D_{1}(\varepsilon)$ and $D_{2}(\varepsilon)$ don't intersect
$L_{t_{3}}$, so $L_{t_{3}}$ must self intersect. This implies
that there is no CMC foliation on $\H^{3}/\Gamma$. 
The claim follows.

Therefore, there exists a quasi-Fuchsian $3$-manifold which does 
not admit CMC foliations.

\bibliographystyle{amsalpha}
\bibliography{foliationbio}

\end{document}